\documentclass{amsart}
\usepackage{amsfonts,latexsym,amsmath,amscd,amssymb}
\usepackage{bbm}

\theoremstyle{plain}
\newtheorem{theorem}{Theorem}
\numberwithin{theorem}{section}

\newtheorem{corollary}{Corollary}
\numberwithin{corollary}{section}

\newtheorem{definition}{Definition}
\numberwithin{definition}{section}

\newtheorem{lemma}{Lemma}
\numberwithin{lemma}{section}

\numberwithin{proposition}{section}

\newtheorem{remark}{Remark}
\numberwithin{remark}{section}

 \numberwithin{equation}{section}

\newcommand {\be}{\begin{equation}}
\newcommand {\ee}{\end{equation}}

\newcommand{\h}{\begin{eqnarray*}}
 \newcommand{\e}{\end{eqnarray*}}

\begin{document}
\title{Even dimensional manifolds and generalized anomaly cancellation formulas}
\author{ Fei Han}
\address{F. Han, \ Department of Mathematics, University of California,
Berkeley, CA, 94720-3840} \email{feihan@math.berkeley.edu}

\author{Xiaoling Huang}
\address{X. Huang, \  Department of Mathematics, University of California,
Santa Barbara, CA, 93106} \email{xiaoling@math.ucsb.edu}
\date{September 1, 2005}


\maketitle

\begin{abstract} We give a direct proof of a cancellation formula raised in
[7] on the level of differential forms. We also obtain more
cancellation formulas for even dimensional Riemannian manifolds
with a complex line bundle involved. Relations among these
cancellation formulas are discussed.
\end{abstract}

\section {Introduction}
In 1983, the physicists Alvarez-Gaum\'e and Witten [1] discovered
the "miraculous cancellation" formula for gravitational anomaly
which reveals a beautiful relation between the top components of
the Hirzebruch $\widehat{L}$-form and $\widehat{A}$-form of a
12-dimensional smooth Riemannian manifold $M$ as follows,\newline
\be
\left\{\widehat{L}(TM,\nabla^{TM})\right\}^{(12)}=\left\{8\widehat{A}
(TM,\nabla^{TM})\mathrm{ch}(T_\mathbb{C}M,\nabla^{T_\mathbb{C}M})
-32\widehat{A}(TM,\nabla^{TM})\right\}^{(12)},\ee where
$T_\mathbb{C}M$ denotes the complexification of $TM$ and
$\nabla^{T_\mathbb{C}M}$ is canonically induced from
$\nabla^{TM}$, the Levi-Civita connection associated to the
Riemannian structure of $M.$

Kefeng Liu [10] established higher dimensional ``miraculous
cancellation" formulas for $(8k+4)$-dimensional Riemannian
manifolds by developing modular invariance properties of
characteristic forms. In [10], he proved that for each
$(8k+4)$-dimensional smooth Riemannian manifold $M$ the following
identity holds,

\be \Big\{\widehat{L}(TM,\nabla^{TM})\Big\}^{(8k+4)}=8\sum_{j=0}^k
2^{6k-6j}\Big\{\widehat{A}(TM,\nabla^{TM})
\mathrm{ch}b_j\Big\}^{(8k+4)},\ee where the $b_j$'s are elements
in $KO(M)\otimes {\mathbb C}$. Liu's formula refines the argument
of Landweber [9] to the level of differential forms and is a
higher dimensional generalization of (1.1). One can also use (1.2)
to deduce the Ochanine divisibility [12] from the
Atiyah-Hirzebruch divisibility [3] for $(8k+4)$-dimensional smooth
closed spin manifolds. In fact, the Atiyah-Hirzebruch divisibility
guarantees that $
\langle\widehat{A}(TM,\nabla^{TM})\mathrm{ch}(E\otimes\mathbb{C}),
[M]\rangle $ is even when $M$ is a smooth closed
$(8k+4)$-dimensional spin manifold and $E$ is a real vector bundle
on $M$. Thus (1.2) implies the signature of $M$ is divisible by
16. This is just the Ochanine divisibility, which generalizes the
famous Rokhlin divisibility for spin 4-manifolds (when $k=0)$.

In [6, 7], for each $(8k+4)$-dimensional smooth Riemannian
manifold $M$, a more general cancellation formula that involves a
complex line bundle on $M$ is established. To be precise, the
authors proved that for each $(8k+4)$-dimensional smooth
Riemannian manifold $M$ and a complex line bundle $\xi$ on $M$, or
equivalently a rank 2 real oriented bundle on $M$, the following
identity holds,

\be
\left\{\frac{\widehat{L}(TM,\nabla^{TM})}{\cosh^2\left(\frac{e}{2}\right)}\right\}^{(8k+4)}
=8\sum_{j=0}^k 2^{6k-6j}\left\{\widehat{A}(TM,\nabla^{TM})
\mathrm{ch}b_j'\cosh\left(\frac{e}{2}\right)\right\}^{(8k+4)}, \ee
where the $b_j'$s are elements in $KO(M)\otimes \mathbb C$,
dependent on $(M,\nabla^{TM})$ and $(\xi,\nabla^{\xi});
e=e(\xi,\nabla^{\xi})$ is the Euler form of $(\xi,\nabla^{\xi}).$
Putting $k=1$ in(1.3), one has

\be \left\{ {\widehat{L}(TM, \nabla^{TM})\over \cosh^2({c\over
2})}\right\}^{(12)}=\left\{ \left[ 8\widehat{A}(TM,
\nabla^{TM}){\rm ch} ( T_{\mathbbm{C}}M,\nabla^{T_{\mathbbm{C}
}M})-32\widehat{A}(TM, \nabla^{TM})\right.\right.\ee
$$- 24\left.\left.\widehat{A}(TM, \nabla^{TM})\left(e^c+e^{-c}-2\right)\right]
\cosh\left({c\over 2}\right)\right\}^{(12)}.$$ This is a twisted
version of the original miraculous cancellation formula (1.1).
When formula (1.3) is applied to spin$^c$ manifolds, the authors
are led directly to a refined version of [11, Theorem 4.2], which
is a beautiful analytic version of the Ochanine congruence formula
[12].

In [7], to obtain a direct proof of [11, Theorem 4.1], an analytic
version of the Finashin congruence formula [5], the authors
applied the following identity, \be
\begin{split}\frac{1}{8}&\int_B
\widehat{L}(TB,\nabla^{TB})\frac{\sinh\left(\frac{e}{2}\right)}
{\cosh\left(\frac{e}{2}\right)}\\
=&\sum_{r=0}^k 2^{6k-6r}\int_B \widehat{A}(TB,\nabla^{TB})
\Big(\mathrm{ch}(b_r(T_\mathbb{C}B+N_\mathbb{C},
\mathbb{C}^2))\\
&\left.-\cosh\left(\frac{e}{2}\right)\mathrm{ch}(b_r(T_\mathbb{C}B+N_\mathbb{C},
N_\mathbb{C}))\right)\frac{1}{2\sinh\left(\frac{e}{2}\right)},
\end{split}\ee
 where
$(B,\nabla^{TB})$ is an $(8k+2)$-dimensional smooth Riemannian
manifold, $(N,\nabla^{N})$ is a rank two real oriented Euclidean
vector bundle on $B$ and $e=e(N,\nabla^N)$ is the associated Euler
form of $(N,\nabla^N)$. This identity is very crucial in their
proof and they proved it by using the cobordism argument. They
also pointed out that (1.5) can be refined to the level of
differential forms and one should be able to prove this directly
by still using the modular invariance method without passing to
the cobordism argument. One of the purposes of this article is to
refine (1.5) to the level of differential forms (Theorem 3.1) and
give such a direct proof. We also obtain a analogous formula for
$(8k+6)$ dimensional Riemannian manifold (Theorem 3.2). One can
view Theorem 3.1 and Theorem 3.2 as generalized miraculous
cancellation formulas on $(8k+2)$ and $(8k+6)$ dimensional smooth
Riemannian manifolds respectively.

With a twisting complex line bundle, we also obtain a unified
cancellation formula (Theorem 3.3) for each even dimensional
smooth Riemannian manifold via the same argument. When the
manifold is of dimension $(8k+4)$ and the bundle is trivial, our
cancellation formula becomes Liu's cancellation formula (1.2).
This unified formula is still a product of the modular invariance
method developed in [10]. Finally, on the level of characteristic
numbers, we discuss relations among cancellation formulas on
manifolds of different dimensions by applying the method of
integration along the fibre.

\section{Modular invariance and characteristic forms}
The purpose of this section is to review the necessary knowledge
on characteristic forms and modular forms that we are going to
use. We also briefly review cancellation formulas obtained in [10]
and [6, 7].

\subsection{Characteristic forms}
Let $M$ be an even-dimensional smooth Riemanniann manifold. Let
$\nabla^{TM}$ be the associated Levi-Civita connection and
$R^{TM}=(\nabla^{TM})^2$ be the curvature of $\nabla^{TM}$ . Let
$\widehat{A}(TM, \nabla^{TM})$, $\widehat{L}(TM, \nabla^{TM})$ be
the Hirzebruch characteristic forms defined respectively by (cf.
[15])
\begin{equation}
\begin{split}
&\widehat{A}(TM, \nabla^{TM}) ={\det}^{1/2}\left({{\sqrt{-1}\over
4\pi}R^{TM} \over \sinh\left({ \sqrt{-1}\over
4\pi}R^{TM}\right)}\right), \\ &\widehat{L}(TM, \nabla^{TM})
={\det}^{1/2}\left({{\sqrt{-1}\over 2\pi}R^{TM} \over \tanh\left({
\sqrt{-1}\over 4\pi}R^{TM}\right)}\right).
\end{split}
\end{equation}

Let $E$, $F$ be two Hermitian vector bundles over $M$ carrying
Hermitian connections $\nabla^E$, $\nabla^F$ respectively. Let
$R^E=(\nabla^{E})^2$ (resp. $R^F=(\nabla^{F})^2$) be the curvature
of $\nabla^E$ (resp. $\nabla^F$). If we set the formal difference
$G=E-F$, then  $G$ carries an induced Hermitian connection
$\nabla^G$ in an obvious sense. We define the associated Chern
character form  as (cf. [15]) \begin{equation}{\rm
ch}(G,\nabla^G)={\rm tr}\left[\exp\left({\sqrt{-1}\over
2\pi}R^E\right)\right] -{\rm tr}\left[\exp\left({\sqrt{-1}\over
2\pi}R^F\right)\right].\end{equation}

In the rest of this paper, for simplicity, when there are no
confusion about the Hermitian connection $\nabla^E$ on a Hermitian
vector bundle $E$, we will only write ${\rm ch}(E)$ for the
associated Chern character form.

For any complex number $t$, let
$$\Lambda_t(E)={\mathbb C}|_M+tE+t^2\Lambda^2(E)+\cdots ,
\\\ S_t(E)={\mathbb C}|_M+tE+t^2S^2(E)+\cdots$$  denote respectively
the total exterior and symmetric powers  of $E$, which live in
$K(M)[[t]].$ The following relations between these two operations
[2, Chap. 3] hold, \be S_t(E)=\frac{1}{\Lambda_{-t}(E)},\ \ \ \
 \Lambda_t(E-F)=\frac{\Lambda_t(E)}{\Lambda_t(F)}.\ee
Moreover, if $\{\omega_i \}$, $\{{\omega_j}' \}$ are formal Chern
roots for Hermitian vector bundles $E$, $F$ respectively, then [8,
Chap. 1] \be
\mathrm{ch}(\Lambda_t{(E)})=\prod\limits_i(1+e^{\omega_i}t).\ee
Therefore,  we have the following formulas for Chern character
forms, \be{\rm ch}(S_t(E) )=\frac{1}{{\rm ch}(\Lambda_{-t}(E)
)}=\frac{1}{\prod\limits_i (1-e^{\omega_i}t)}\ ,\ee \be{\rm
ch}(\Lambda_t(E-F) )=\frac{{\rm ch}(\Lambda_t(E) )}{{\rm
ch}(\Lambda_t(F)
)}=\frac{\prod\limits_i(1+e^{\omega_i}t)}{\prod\limits_j(1+e^{{\omega_j}'}t)}\
.\ee

If $W$ is a  real Euclidean vector bundle over $M$ carrying a
Euclidean connection $\nabla^W$, its complexification
$W_\mathbb{C}=W\otimes \mathbb{C}$ is a complex vector bundle over
$M$ carrying a canonically induced Hermitian metric from the
Euclidean metric of $W$ as well as a Hermitian connection
$\nabla^{W_\mathbb{C}}$ induced from $\nabla^W$.

\subsection{Results needed on the Jacobi theta functions and
modular forms} The four Jacobi theta functions are defined as
follows (cf. [4]): \be\theta(v,\tau)=2q^{1/8}\sin(\pi v)
\prod_{j=1}^\infty\left[(1-q^j)(1-e^{2\pi
\sqrt{-1}v}q^j)(1-e^{-2\pi \sqrt{-1}v}q^j)\right]\ ,\ee \be
\theta_1(v,\tau)=2q^{1/8}\cos(\pi v)
 \prod_{j=1}^\infty\left[(1-q^j)(1+e^{2\pi \sqrt{-1}v}q^j)
 (1+e^{-2\pi \sqrt{-1}v}q^j)\right]\ ,\ee
\be \theta_2(v,\tau)=\prod_{j=1}^\infty\left[(1-q^j)
 (1-e^{2\pi \sqrt{-1}v}q^{j-1/2})(1-e^{-2\pi \sqrt{-1}v}q^{j-1/2})\right]\
 ,\ee
\be \theta_3(v,\tau)=\prod_{j=1}^\infty\left[(1-q^j) (1+e^{2\pi
\sqrt{-1}v}q^{j-1/2})(1+e^{-2\pi \sqrt{-1}v}q^{j-1/2})\right]\
,\ee where $q=e^{2\pi \sqrt{-1}\tau}$ with $\tau\in {\mathbb H}$,
the upper half complex plane.

Let \be {\theta}'(0,\tau)=\left. {\partial \theta (v,\tau)\over
\partial v}\right|_{v=0}.\ee Then the following Jacobi
identity (cf. [4]) holds, \be {\theta}'(0,\tau)= \pi\,
\theta_1(0,\tau)\theta_2(0,\tau)\theta_3(0,\tau)\, .\ee

Denote $SL_2(\mathbb{Z})= \left\{\left. \left(\begin{array}{cc}a&b\\
c&d\end{array}\right)\, \right|\ a,b,c,d\in\mathbb{Z},\ ad-bc=1\right\}$ the modular group. Let $S=\left(\begin{array}{cc}0&-1\\
1&0\end{array}\right)$, $T=\left(\begin{array}{cc}1&1\\
0&1\end{array}\right)$ be the two generators of
$SL_2(\mathbb{Z})$. They act on ${\mathbb H}$ by $S\tau =-1/\tau$,
$T\tau=\tau+1$. One has the following transformation laws of theta
functions under the actions of $S$ and $T$ (cf. [4]): \be
\theta(v,\tau+1)=e^{\pi \sqrt{-1}\over 4}\theta(v,\tau),\ \ \
\theta\left(v,-{1}/{\tau}\right)={1\over\sqrt{-1}}\left({\tau\over
\sqrt{-1}}\right)^{1/2} e^{\pi\sqrt{-1}\tau v^2}\theta\left(\tau
v,\tau\right)\ ;\ee \be \theta_1(v,\tau+1)=e^{\pi \sqrt{-1}\over
4}\theta_1(v,\tau),\ \ \
\theta_1\left(v,-{1}/{\tau}\right)=\left({\tau\over
\sqrt{-1}}\right)^{1/2} e^{\pi\sqrt{-1}\tau v^2}\theta_2(\tau
v,\tau)\ ;\ee \be\theta_2(v,\tau+1)=\theta_3(v,\tau),\ \ \
\theta_2\left(v,-{1}/{\tau}\right)=\left({\tau\over
\sqrt{-1}}\right)^{1/2} e^{\pi\sqrt{-1}\tau v^2}\theta_1(\tau
v,\tau)\ ;\ee \be\theta_3(v,\tau+1)=\theta_2(v,\tau),\ \ \
\theta_3\left(v,-{1}/{\tau}\right)=\left({\tau\over
\sqrt{-1}}\right)^{1/2} e^{\pi\sqrt{-1}\tau v^2}\theta_3(\tau
v,\tau)\ ;\ee
\be \theta'(0,\tau+1)=e^{\pi \sqrt{-1}\over 4}\theta'(0,\tau),\\\
\theta'(0,-\frac{1}{\tau})=\frac{1}{\sqrt{-1}}\tau^{3/2}\theta'(0,\tau).\ee

\begin{definition}A modular form over $\Gamma$, a subgroup of
$SL_2(\mathbb{Z})$, is a holomorphic function $f(\tau)$ on
${\mathbb H}\cup\{\infty\}$ such that \be f(g\tau)
:=f\left(\frac{a\tau+b}{c\tau+d}\right)
=\chi(g)(c\tau+d)^kf(\tau),\ \ \ \  \forall \
 g=\left(\begin{array}{cc}a&b\\
c&d\end{array}\right)\in \Gamma,\ee
 where $\chi:\Gamma\rightarrow\mathbb{C}^*$ is a character of
 $\Gamma$. $k$ is called the weight of $f$.\end{definition}

Denote by $\theta_j=\theta_j(0,\tau)$, $1\leq j\leq 3$, and define
\be \delta_1(\tau)={1\over 8}(\theta_2^4+\theta_3^4),\ \ \ \
\varepsilon_1(\tau)={1\over 16}\theta_2^4\theta_3^4,\ee
\be\delta_2(\tau)=-{1\over 8}(\theta_1^4+\theta_3^4),\ \ \ \
\varepsilon_2(\tau)={1\over 16}\theta_1^4\theta_3^4.\ee
 They admit Fourier expansion (cf. [9])
\be \delta_1(\tau)={1\over 4}+6q+\cdots,\ \ \ \
\varepsilon_1(\tau)={1\over 16}-q+\cdots,\ee \be
\delta_2(\tau)=-{1\over 8}-3q^{1/2}+\cdots,\ \ \ \
\varepsilon_2(\tau)=q^{1/2}+\cdots,\ee where the ``$\cdots$" terms
are higher degree terms all having integral coefficients. They
also satisfy the following transformation laws under $S$ (cf. [9]
and [10]), \be\delta_2(-{1/ \tau}) =\tau^2\delta_1(\tau),\ \ \ \
\varepsilon_2(-{1/ \tau})=\tau^4\varepsilon_1(\tau).\ee

Let $\Gamma_0(2)$, $\Gamma^0(2)$ be the two subgroups of
$SL_2(\mathbb{Z})$ defined by
$$ \Gamma_0(2)=\left\{\left.\left(\begin{array}{cc}
a&b\\
c&d
\end{array}\right)\in SL_2(\mathbb{Z})\,\right|\,c\equiv 0\ \ {\rm mod} \ \ 2{\mathbb Z}\right\},$$

$$ \Gamma^0(2)=\left\{\left.\left(\begin{array}{cc}
a&b\\
c&d
\end{array}\right)\in SL_2(\mathbb{Z})\,\right|\,b\equiv 0\ \ {\rm mod} \ \ 2{\mathbb Z}\right\}.$$
Then $T,\ ST^2ST$ are the two generators of $\Gamma_0(2)$, while
$STS,\ T^2STS$ are the two generators of $\Gamma^0(2)$.

The following weaker version of [10, Lemma 2] will be used in the
next section.

\begin{lemma} One has that $\delta_2$ (resp.
$\varepsilon_2$) is a modular form of weight $2$ (resp. $4$) over
$\Gamma^0(2)$. Furthermore, ${\mathcal M}_{\mathbb R}(\Gamma^0(2)
)={\mathbb R}[\delta_2(\tau),\varepsilon_2(\tau)]$, where
${\mathcal M}_{\mathbb R}(\Gamma)$ denote the ring of modular
forms over $\Gamma$ with real Fourier coefficients.\end{lemma}

\subsection{Cancellation formulas for ${(8k+4)}$-dimensional
Riemannian manifolds} Let $M$ be an ${(8k+4)}$-dimensional smooth
Riemannian manifold with Levi-Civita connection $\nabla^{TM}$. Let
$\nabla^{T_{\mathbb C} M}$ be the canonically induced Hermitian
connection on $T_\mathbb{C}M=TM\otimes\mathbb{C}$. Let $V$ be a
rank $2l$ Euclidean vector bundle over $M$ carrying a Euclidean
connection $\nabla^{V}$. Let $\xi$ be a rank two oriented
Euclidean vector bundle carrying a Euclidean connection
$\nabla^{\xi}$. Let $\nabla^{\xi_\mathbb{C}}$ be the canonically
induced Hermitian connection on $\xi_\mathbb{C}=\xi\otimes
\mathbb{C}$. Let $c=e(\xi, \nabla^{\xi})$ be the Euler form of
$\xi$ canonically associated to $\nabla^{\xi}$. If $W$ is a
complex vector bundle over $M$, denote
$\widetilde{W}=W-\mathbb{C}^{{\rm dim}_\mathbb{C}W}|_M$ in $K(M)$.

Using the same notations as in Section 2.1, we construct two
formal power series in $q^{1/2}$ with coefficients in the
semi-group generated by complex vector bundles over $M$, which are
introduced in [6, 7] to prove Theorem 2.1 in this text, \be
\begin{split}\Theta_1(T_{\mathbb C}M,V_{\mathbb C},\xi_{\mathbb
C})&=\bigotimes_{n=1}^\infty S_{q^n}(\widetilde{T_{\mathbb C}M})
\otimes \bigotimes_{m=1}^\infty
\Lambda_{q^m}(\widetilde{V}_{\mathbb C}-2\widetilde{\xi}_{\mathbb
C})\\
&\otimes \bigotimes_{r=1}^\infty\Lambda_{q^{r-{1\over
2}}}(\widetilde{\xi}_{\mathbb
C})\otimes\bigotimes_{s=1}^\infty\Lambda_{-q^{s-{1\over
2}}}(\widetilde{\xi}_{\mathbb C}),\end{split}\ee

\be \begin{split}\Theta_2(T_{\mathbb C}M,V_{\mathbb
C},\xi_{\mathbb C})&=\bigotimes_{n=1}^\infty
S_{q^n}(\widetilde{T_{\mathbb C}M}) \otimes
\bigotimes_{m=1}^\infty \Lambda_{-q^{m-{1\over
2}}}(\widetilde{V}_{\mathbb C}-2\widetilde{\xi}_{\mathbb C})\\
&\otimes \bigotimes_{r=1}^\infty\Lambda_{q^{r-{1\over
2}}}(\widetilde{\xi}_{\mathbb
C})\otimes\bigotimes_{s=1}^\infty\Lambda_{q^{s}}(\widetilde{\xi}_{\mathbb
C}).\end{split}\ee $\Theta_1(T_{\mathbb C}M,V_{\mathbb
C},\xi_{\mathbb C })$ and $\Theta_2(T_{\mathbb C}M,V_{\mathbb
C},\xi_{\mathbb C})$ admit formal Fourier expansion  in $q^{1/2}$
as

\be\Theta_1(T_{\mathbb C}M,V_{\mathbb C},\xi_{\mathbb
C})=A_0(T_{\mathbb C}M,V_{\mathbb C},\xi_{\mathbb C})+
A_1(T_{\mathbb C}M,V_{\mathbb C},\xi_{\mathbb
C})q^{1/2}+\cdots,\ee

\be \Theta_2(T_{\mathbb C}M,V_{\mathbb C},\xi_{\mathbb
C})=B_0(T_{\mathbb C}M,V_{\mathbb C},\xi_{\mathbb C})
+B_1(T_{\mathbb C}M,V_{\mathbb C},\xi_{\mathbb
C})q^{1/2}+\cdots,\ee where the $A_j$'s and $B_j$'s are elements
in the semi-group formally generated by Hermitian vector bundles
over $M$. Moreover, they carry canonically induced Hermitian
connections denoted by $\nabla^{A_j}$ and $\nabla^{B_j}$
respectively, and $\nabla^{\Theta_i(M, V, \xi)}$ are the induced
Hermitian connections with $q^{1/2}$-coefficients on $\Theta_i$
from the $\nabla^{A_j}$ and $\nabla^{B_j}$.

Now, we can state a cancellation formula, which is obtained in [6,
7].

\begin{theorem}[Han-Zhang, 2003] If the equality for the
first Pontrjagin forms $p_1(TM,\nabla^{TM})=p_1(V,\nabla^V)$
holds, then one has an equality for $(8k+4) $-forms, \be \left\{
{\widehat{A}(TM,\nabla^{TM}){\det}^{1/2}
\left(2\cosh\left({\sqrt{-1}\over 4\pi}R^V\right)\right) \over
\cosh^2({e\over 2})}\right\}^{(8k+4)} \ee
$$=2^{l+2k+1}\sum_{r=0}^k 2^{-6r}\left\{
\widehat{A}(TM,\nabla^{TM}){\rm ch}( b_r(T_{\mathbb C}M,V_{\mathbb
C},\xi_{\mathbb C }) )\cosh\left({e\over
2}\right)\right\}^{(8k+4)}, $$ where each $b_r(T_\mathbb{C} M,
V_\mathbb{C}, \xi_\mathbb{C})$, $0\leq r\leq k$, is a  canonical
integral linear combination of $B_j(T_{\mathbb C}M, V_{\mathbb C},
\xi_{\mathbb C }) $, $0\leq j\leq r$.\end{theorem}

From now on, denote $\Theta_i(T_\mathbb{C} M, T_\mathbb{C} M,
\xi_\mathbb{C})$ as $\Theta_i(T_\mathbb{C} M, \xi_\mathbb{C})$ and
$b_r(T_\mathbb{C} M, T_\mathbb{C} M, \xi_\mathbb{C})$ as
$b_r(T_\mathbb{C} M, \xi_\mathbb{C})$.

Taking $V=TM$ in (2.28), we have

\begin{corollary} The following identity of characteristic
forms holds, \be \left\{\frac{\widehat{L}(TM,\nabla^{TM})}
{\cosh^2\left(\frac{e}{2}\right)}\right\}^{(8k+4)}\ee
$$=8\sum_{r=0}^k 2^{6k-6r}\left\{\widehat{A}(TM,\nabla^{TM}){\rm
ch}(b_r(T_\mathbb{C}M,
\xi_\mathbb{C}))\cosh\left(\frac{e}{2}\right)\right\}^{(8k+4)}.$$
\end{corollary}
Moreover, taking $\xi =\mathbb{R}^2$ in (2.29), we have

\begin{corollary} The following identity of characteristic
forms holds, \be
\left\{\widehat{L}(TM,\nabla^{TM})\right\}^{(8k+4)}\ee
$$=8\sum_{r=0}^k 2^{6k-6r}\left\{\widehat{A}(TM,\nabla^{TM}){\rm
ch}(b_r(T_\mathbb{C}M, \mathbb{C}^2))\right\}^{(8k+4)}.$$
\end{corollary}
(2.30) is exactly the cancellation formula obtained in [10].

\section{Cancellation formulas for even dimensional Riemannian manifolds}
Let $B$ be an $(8k+2)$-dimensional smooth oriented Riemannian
manifold with Levi-Civita connection $\nabla^{TB}$. let $\pi: N\to
B $ be a rank two real oriented Euclidean vector bundle over $B$
carrying the Euclidean connection $\nabla^N$. Let
$R^N=(\nabla^{N})^2$ be the curvature of $\nabla^N$ and
$e=e(N,\nabla^{N})$ be the Euler form of $(N,\nabla^N)$. Then we
have the following cancellation formula for characteristic numbers
\be
\begin{split}\frac{1}{8}&\int_B
\widehat{L}(TB,\nabla^{TB})\frac{\sinh\left(\frac{e}{2}\right)}
{\cosh\left(\frac{e}{2}\right)}\\
=&\sum_{r=0}^k 2^{6k-6r}\int_B \widehat{A}(TB,\nabla^{TB}) \Big
(\mathrm{ch} (b_r(T_\mathbb{C}B+N_\mathbb{C},
\mathbb{C}^2))\\
&\left.-\cosh\left(\frac{e}{2}\right)\mathrm{ch}(b_r(T_\mathbb{C}B+N_\mathbb{C},
N_\mathbb{C}))\right)\frac{1}{2\sinh\left(\frac{e}{2}\right)},
\end{split}\ee which is applied in [7] to give the analytic
Finashin congruence [11] a direct analytic proof via a beautiful
Rokhlin type congruence formula obtained in [14]. In fact, formula
(3.1) holds on the level of differential forms. In \S 3.1 we prove
the form-level version of (3.1) directly by applying the modular
invariance argument. Also in this subsection, we give the
$(8k+6)$-dimensional analogue without proof.

In \S 3.2, for all even dimensional smooth Riemannian manifolds,
we obtain a general type of cancellation formulas, which imply
Liu's formula (2.30) as a special case.

\subsection{The direct proof for the form-level version of (3.1)}
We make the same assumptions and use the same notations as in \S
2.3. Define \be \begin{split}
\Theta_{1}(T_{\mathbb{C}}B+N_{\mathbb{C}},N_{\mathbb{C}})
&=\bigotimes_{n=1}^\infty\,S_{q^n}(\widetilde{T_{\mathbb{C}}B+N_{\mathbb{C}}})
 \otimes\bigotimes_{m=1}^\infty\Lambda_{q^m}(\widetilde{T_{\mathbb{C}}B+N_{\mathbb{C}}}
-2\widetilde{N_{\mathbb{C}}})\\
 &\otimes\bigotimes_{r=1}^\infty\Lambda_{q^{(r-\frac{1}{2})}}(\widetilde{N_{\mathbb{C}}})
 \otimes\bigotimes_{s=1}^\infty\Lambda_{-q^{(s-\frac{1}{2})}}(\widetilde{N_{\mathbb{C}}}),\end{split} \ee
\be
\begin{split}\Theta_{2}(T_{\mathbb{C}}B+N_{\mathbb{C}},N_{\mathbb{C}})
&=\bigotimes_{n=1}^\infty\,S_{q^n}(\widetilde{T_{\mathbb{C}}B+N_{\mathbb{C}}})
 \otimes\bigotimes_{m=1}^\infty\Lambda_{-q^{(m-\frac{1}{2})}}(\widetilde{T_{\mathbb{C}}B+N_{\mathbb{C}}}
-2\widetilde{N_{\mathbb{C}}})\\
 &\otimes\bigotimes_{r=1}^\infty\Lambda_{q^{(r-\frac{1}{2})}}(\widetilde{N_{\mathbb{C}}})
 \otimes\bigotimes_{s=1}^\infty\Lambda_{q^s}(\widetilde{N_{\mathbb{C}}}). \end{split}
 \ee
and denote $A_i(T_{\mathbb{C}}B+N_{\mathbb{C}},N_{\mathbb{C}})$
and $B_i(T_{\mathbb{C}}B+N_{\mathbb{C}},N_{\mathbb{C}})$ the
coefficients in their Fourier expansions respectively. Then we
have the following equality for ${(8k+2)}$-forms associated to $B$
and $N$.

\begin{theorem}The following identity holds,
\be
\left\{\widehat{L}(TB,\nabla^{TB})\frac{\sinh\left(\frac{e}{2}\right)}
{\cosh\left(\frac{e}{2}\right)}\right\}^{(8k+2)} =8\sum_{r=0}^k
2^{6k-6r}h_r,\ee where each $h_r$, $0\leq r\leq k$, is a canonical
integral linear combination of the characteristic forms
$$
\left\{\widehat{A}(TB,\nabla^{TB})\frac{
\mathrm{ch}\left(B_j(T_\mathbb{C}B+N_\mathbb{C},\mathbb{C}^2)\right)-\cosh\left(\frac
e2\right)\mathrm{ch}\left(B_j(T_\mathbb{C}B+N_\mathbb{C},
N_\mathbb{C})\right)
}{2\sinh\left(\frac{e}{2}\right)}\right\}^{(8k+2)}, $$ $0\leq
j\leq r.$ Actually, $h_r$ is just
$$\left\{\widehat{A}(TB,\nabla^{TB}){{\rm ch}\left( b_r(T_{\mathbb
C}B+N_{\mathbb C},{\mathbb C }^2)\right) -\cosh\left({e\over
2}\right){\rm ch}\left( b_r(T_{\mathbb C}B+N_{\mathbb
C},N_{\mathbb C })\right)\over 2\sinh\left({e\over
2}\right)}\right\}^{(8k+2)}.$$ \end{theorem}

\noindent{\it Proof}. As in [10], we use the formal Chern roots
$\{\pm 2\pi \sqrt{-1}x_j\}$ for
$(T_{\mathbb{C}}B,\nabla^{T_{\mathbb{C}}B})$ Let
$e=2\pi\sqrt{-1}u, q=e^{2\pi\sqrt{-1}\tau}$ with $\tau\in
\mathbb{H}$, the upper half complex plane. Set \be \begin{split}
{Q}_{1}(\tau)=&\widehat{L}(TB,\nabla^{TB})
\frac{\cosh\left(\frac{e}{2}\right)}{\sinh\left(\frac{e}{2}\right)}\\
&
\cdot\left(\mathrm{ch}\left(\Theta_{1}(T_{\mathbb{C}}B+N_{\mathbb{C}},{\mathbb{C}}^2)\right)
-\frac{\mathrm{ch}\left(\Theta_{1}(T_{\mathbb{C}}B+N_{\mathbb{C}},N_{\mathbb{C}})\right)}
{\cosh^2\left(\frac{e}{2}\right)}\right)\\
=&2^{4k+1}\left(\prod_{j=1}^{4k+1}\frac{\pi x_j}{\sin (\pi
x_j)}\right)\left(\prod_{j=1}^{4k+1}\cos(\pi x_j)\right)
\frac{\cos(\pi u)}{\sin(\pi u)}\\
&\cdot
\left(\mathrm{ch}\left(\Theta_{1}(T_{\mathbb{C}}B+N_{\mathbb{C}},{\mathbb{C}}^2)\right)
-\frac{\mathrm{ch}\left(\Theta_{1}(T_{\mathbb{C}}B+N_{\mathbb{C}},N_{\mathbb{C}})\right)}
{\cos^2\left(\pi u\right)}\right), \end{split}\ee

\be \begin{split}{Q}_{2}(\tau)=&\widehat{A}(TB,\nabla^{TB})
\frac{1}{2\sinh\left(\frac{e}{2}\right)}
\Big(\mathrm{ch}\left(\Theta_{2}(T_{\mathbb{C}}B+N_{\mathbb{C}},{\mathbb{C}}^2)\right)\\
&\left.-\cosh\left(\frac{e}{2}\right)
\mathrm{ch}\left(\Theta_2(T_\mathbb{C}B+N_\mathbb{C},N_\mathbb{C}\right)\right)\\
=& \left(\prod_{j=1}^{4k+1}\frac{\pi x_j}{\sin (\pi x_j)}\right)
\frac{1}{2\sin (\pi u)}
\left(\mathrm{ch}\left(\Theta_{2}(T_{\mathbb{C}}B+N_{\mathbb{C}},{\mathbb{C}}^2\right)\right)\\
&-\cos (\pi u)
\mathrm{ch}\left(\Theta_2(T_\mathbb{C}B+N_\mathbb{C},N_\mathbb{C}\right)))
.\end{split}\ee

We can actually write $Q_1(\tau)$ and $Q_2(\tau)$ in terms of the
Jacobi theta-functions as \be \begin{split}
Q_1(\tau)=&2^{4k+1}\left(\prod_{j=1}^{4k+1}x_j\frac{\theta'(0,\tau)}{\theta(x_j,\tau)}
\frac{\theta_{1}(x_j,\tau)}{\theta_{1}(0,\tau)}\right)
\frac{\theta'(0,\tau)}{\theta(u,\tau)}\\
&\cdot\left(\frac{\theta_1(u,\tau)}{\theta_1(0,\tau)}
-\frac{\theta_{1}(0,\tau)}{\theta_{1}(u,\tau)}
\frac{\theta_{3}(u,\tau)}{\theta_{3}(0,\tau)}
\frac{\theta_{2}(u,\tau)}{\theta_{2}(0,\tau)}\right)\end{split}
\ee and

\be \begin{split} Q_2(\tau)=&\frac 12\left(\prod_{j=1}^{4k+1}x_j
\frac{\theta'(0,\tau)}{\theta(x_j,\tau)}
\frac{\theta_2(x_j,\tau)}{\theta_2(0,\tau)}\right)
\frac{\theta'(0,\tau)}{\theta(u,\tau)}\\
&\cdot\left(\frac{\theta_2(u,\tau)}{\theta_2(0,\tau)}
-\frac{\theta_2(0,\tau)}{\theta_2(u,\tau)}
\frac{\theta_3(u,\tau)}{\theta_3(0,\tau)}
\frac{\theta_1(u,\tau)}{\theta_1(0,\tau)}\right)\end{split}\ee In
fact, by (3.2), (2.5) and (2.6), \be
\mathrm{ch}\left(\Theta_{1}(T_\mathbb{C}B+N_\mathbb{C},N_\mathbb{C})\right)=
\prod\limits_{n=1}^\infty\frac{\mathrm{ch}\Lambda_-q^n(\mathbb{C}^{8k+2})}{\mathrm{ch}\Lambda_-q^n(T_\mathbb{C}B)}
\prod\limits_{n=1}^\infty\frac{\mathrm{ch}\Lambda_-q^n(\mathbb{C}^{2})}{\mathrm{ch}\Lambda_-q^n(N_{\mathbb{C}})}\ee
$$\cdot\prod\limits_{m=1}^\infty\frac{\mathrm{ch}\Lambda_{q^m}(T_{\mathbb{C}}B)}{\mathrm{ch}\Lambda_{q^m}(\mathbb{C}^{8k+2})}
\prod\limits_{m=1}^\infty\frac{\mathrm{ch}\Lambda_{q^m}(\mathbb{C}^2)}{\mathrm{ch}\Lambda_{q^m}(N_{\mathbb{C}})}
\prod\limits_{r=1}^\infty\frac{\mathrm{ch}\Lambda_{q^{(r-\frac{1}{2})}}(N_{\mathbb{C}})}{\mathrm{ch}\Lambda_{q^{(r-\frac{1}{2})}}
(\mathbb{C}^{2})}
\prod\limits_{s=1}^\infty\frac{\mathrm{ch}\Lambda_{-q^{(s-\frac{1}{2})}}(N_\mathbb{C})}{\mathrm{ch}\Lambda_{-q^{(s-\frac{1}{2})}}
(\mathbb{C}^{2})}.$$ From (2.4), the Jacobi identity (2.12) and
(2.7), one deduces directly that \be \prod_{j=1}^{4k+1}\frac{\pi
x_j}{\sin(\pi x_j)}
\prod\limits_{n=1}^\infty\frac{\mathrm{ch}\Lambda_{-q^n}\left(\mathbb{C}^{8k+2}\right)}{\mathrm{ch}\Lambda_{-q^n}(T_{\mathbb{C}}B)}
 \ee $$=\prod_{j=1}^{4k+1}x_j\frac{\pi \theta_1(0,\tau)\theta_2(0,\tau)\theta_3(0,\tau)}{\theta(x_j,\tau)}
 =\prod\limits_{j=1}^{4k+1}x_j\frac{\theta'(0,\tau)}{\theta(x_j,\tau)}.$$

Similarly, from (2.4) and (2.7)-(2.10), one deduces that \be
\prod_{j=1}^{4k+1}\cos(\pi x_j)
\prod\limits_{m=1}^\infty\frac{\mathrm{ch}\Lambda_{q^m}(T_{\mathbb{C}}B)}{\mathrm{ch}\Lambda_{q^m}(\mathbb{C}^{8k+2})}
= \prod_{j=1}^{4k+1}\frac{\theta_1(x_j,\tau)}{\theta_1(0,\tau)},\
\prod\limits_{r=1}^\infty\frac{\mathrm{ch}\Lambda_{q^{(r-\frac{1}{2}})}(N_{\mathbb{C}})}{\mathrm{ch}\Lambda_{q^{(r-\frac{1}{2})}}
(\mathbb{C}^{2})}=\frac{\theta_3(u,\tau)}{\theta(0,\tau)},\ee
$$\frac{1}{\sin{(\pi u})}
\prod\limits_{n=1}^\infty\frac{\mathrm{ch}\Lambda_{-q^n}(\mathbb{C}^{2})}
{\mathrm{ch}\Lambda_{-q^n}(N_\mathbb{C})}
=\frac{\theta'(0,\tau)}{\theta(u,\tau)},\   \frac{1}{\cos{(\pi u})}
\prod\limits_{m=1}^\infty\frac{\mathrm{ch}\Lambda_{q^
m}(\mathbb{C}^{2})}
{\mathrm{ch}\Lambda_{q^m}(N_\mathbb{C})}
=\frac{\theta_1(0,\tau)}{\theta_1(u,\tau)},
$$
$$
\prod\limits_{r=1}^\infty\frac{\mathrm{ch}\Lambda_{-q^
{(r-1/2)}}(N_\mathbb{C})}
{\mathrm{ch}\Lambda_{-q^{(r-1/2)}}(\mathbb{C}^{2})}
=\frac{\theta_1(u,\tau)}{\theta_1(0,\tau)}. $$ Putting (3.5) and
(3.9)-(3.11) together, we get(3.7). By doing similar computations,
one also gets (3.8).

Let $P_1(\tau)=\{Q_1(\tau)\}^{(8k+2)},
P_2(\tau)=\{Q_2(\tau)\}^{(8k+2)}$ be the $(8k+2)$-components of
$Q_1(\tau), Q_2(\tau)$ respectively. Applying the transformation
laws (2.13)-(2.17) to $P_1(\tau)$ and $P_2(\tau)$, we find
$P_1(\tau)$ is a modular form of weight $4k+2$ over $\Gamma_0(2)$;
while $P_2(\tau)$ is a modular form of weight $4k+2$ over
$\Gamma^0(2)$. Moreover, the following identity holds, \be
P_1(-1/\tau)=2^{4k+2}P_2(\tau).\ee

Observe that at any point $x\in M$, up to the volume from
determined by the metric on $T_xM$, both $P_i(\tau), i=1,\ 2$, can
be viewed as a power series of $q^{1/2}$ with real Fourier
coefficients. Thus, one can apply Lemma 2.1 to $P_2(\tau)$ to get,
at $x$, that \be P_2(\tau)
 =h_0(8\delta_2)^{2k+1}+h_0(8\delta_2)^{2k-1}\varepsilon_2
+\cdots+h_k(8\delta_2)\varepsilon_2^k ,\ee where each $h_r$,
$0\leq r\leq k $, is a (canonically) finite integral linear
combination of the forms
$$\left\{\widehat{A}(TB,\nabla^{TB})\frac{\mathrm{ch}\left(B_j(T_\mathbb{C}B+N_\mathbb{C},\mathbb{C}^2)\right)-\cosh\left(\frac
e2\right)\mathrm{ch}\left(B_j(T_\mathbb{C}B+N_\mathbb{C},
N_\mathbb{C})\right)}{2\sinh(\frac{e}{2})} \right\}^{(8k+2)}, $$
$0\leq j\leq r .$

By (2.23) and (3.12), we have \be
\begin{split}P_1(\tau)=&2^{4k+2}\frac{1}{\tau^{4k+2}}P_2(-1/\tau)\\
=&2^{4k+2}\frac{1}{\tau^{4k+2}}\Big[h_0\big(8\delta_2(-1/\tau)\big)^{2k+1}
+h_1\big(8\delta_2(-1/\tau)\big)^{2k-1}\varepsilon_2(-1/\tau)
+\cdots\\
&+h_k\big(8\delta_2(-1/\tau)\big)\big(\varepsilon_2(-1/\tau)\big)^k\Big]\\
=&2^{4k+2}\left[h_0(8\delta_1)^{2k+1}+h_1(8\delta_1)^{2k-1}\varepsilon_1
+\cdots+h_k(8\delta_1)\varepsilon_1^k\right]. \end{split}\ee By
(2.21) and (3.4) and by setting $q=0$ in (3.14), we have
\be\left\{\widehat{L}(TB,\nabla^{TB})\frac{\cosh\left(\frac{e}{2}\right)}
{\sinh\left(\frac {e}{2}\right)}\left(1-\frac{1}{\cosh^2(\frac
{e}{2})}\right)\right\}^{(8k+2)} =2^{6k+3}\sum\limits_{r=0}^k
2^{-6r}h_r. \ee Therefore,
\be
\left\{\widehat{L}(TB,\nabla^{TB})\frac{\sinh\left(\frac{e}{2}\right)}
{\cosh\left(\frac{e}{2}\right)}\right\}^{(8k+2)}
=8\sum_{r=0}^{k}2^{6k-6r}h_r.\ee

We also need to show that each $h_r, 0\leq r\leq k,$ can be
expressed through a canonical integral linear combination of
\be\left\{\widehat{A}(TB,\nabla^{TB})
\frac{1}{2\sinh{\left(\frac{e}{2}\right)}}
\left(\mathrm{ch}\left(B_j(T_\mathbb{C}B+N_\mathbb{C},\mathbb{C}^2)\right)\right.\right.\ee
$$  \left.-\cosh\left(\frac
e2\right)\mathrm{ch}\left(B_j(T_\mathbb{C}B+N_\mathbb{C},
N_\mathbb{C})\right) \right)\Biggr\}^{(8k+2)},$$ $\newline 0\leq
j\leq r,$ with coefficients not
 depending on $x\in M$. As in [10], one can use the induction method to prove this fact
easily by comparing the coefficients of $q^{j/2}, j\geq 0$,
between the two sides of (3.13). For the consideration of the
volumn of this paper, we do not give details here
 but only write down the explicit expressions for $h_0$ and $h_1$ as follows.
\be h_0 =-\left\{\widehat{A}(TB,\nabla^{TB})\frac{1}{2\sinh(\frac
e2)}\left(1-\cosh(\frac e2)\right)\right\}^{(8k+2)}, \ee

\be
h_1=-\left\{\widehat{A}(TB,\nabla^{TB})\frac{1}{2\sinh\left(\frac
e2\right)}\Big[\mathrm{ch}\left(B_1(T_\mathbb{C}B+N_\mathbb{C},\mathbb{C}^2)\right)\right.\ee
$$\left.-\cosh\left(\frac
e2\right)\mathrm{ch}\left(B_1(T_\mathbb{C}B+N_\mathbb{C},N_\mathbb{C})\right)
-24(2k+1)\left(1-\cosh\left(\frac e2
\right)\right)\right]\Biggr\}^{(8k+2)}.
$$

In a summary, we get \be
\left\{\widehat{L}(TB,\nabla^{TB})\frac{\sinh\left(\frac{e}{2}\right)}
{\cosh\left(\frac{e}{2}\right)}\right\}^{(8k+2)} =8\sum_{r=0}^k
2^{6k-6r}h_r, \ee where each $h_r$, $0\leq r\leq k$, is a
canonical integral linear combination of the characteristic forms
$$
\left\{\widehat{A}(TB,\nabla^{TB})\frac{
\mathrm{ch}\left(B_j(T_\mathbb{C}B+N_\mathbb{C},\mathbb{C}^2)\right)-\cosh\left(\frac
e2\right)\mathrm{ch}\left(B_j(T_\mathbb{C}B+N_\mathbb{C},
N_\mathbb{C})\right)
}{2\sinh\left(\frac{e}{2}\right)}\right\}^{(8k+2)}, $$ $0\leq
j\leq r.$

Since both $h_r$'s and $b_r$'s are canonically determined by
induction, one easily finds that,
$h_r=\left\{\widehat{A}(TB,\nabla^{TB}){{\rm ch}\left(
b_r(T_{\mathbb C}B+N_{\mathbb C},{\mathbb C }^2)\right)
-\cosh\left({e\over 2}\right){\rm ch}\left( b_r(T_{\mathbb
C}B+N_{\mathbb C},N_{\mathbb C })\right)\over 2\sinh\left({e\over
2}\right)}\right\}^{(8k+2)}. $
 $\square$

\begin{remark} It's not hard to see that each
$$
\widehat{A}(TB,\nabla^{TB})\frac{
\mathrm{ch}\left(B_j(T_\mathbb{C}B+N_\mathbb{C},\mathbb{C}^2)\right)-\cosh\left(\frac
e2\right)\mathrm{ch}\left(B_j(T_\mathbb{C}B+N_\mathbb{C},
N_\mathbb{C})\right) }{2\sinh\left(\frac{e}{2}\right)}$$  makes
sense as a differential form. \end{remark}

For $(8k+6)$-dimensional manifolds, we have an analogue of Theorem
3.1. Let $B$ be an $(8k+6)$-dimensional smooth oriented Riemannian
manifold and all of the notations in the following theorem make
the same senses as above. Then we can prove verbatim to get

\begin{theorem}The following cancellation formula holds:
$$\left\{\widehat{L}(TB,\nabla^{TB})\frac{\sinh\left(\frac{e}{2}\right)}{\cosh
\left(\frac{e}{2}\right)}\right\}^{(8k+6)}=64\sum_{r=0}^k
2^{6k-6r}h_r ,$$ where each $h_r$, $0\leq r\leq k$, is a canonical
integral linear combination of the characteristic forms
$$
\left\{\widehat{A}(TB,\nabla^{TB})\frac{
\mathrm{ch}\left(B_j(T_\mathbb{C}B+N_\mathbb{C},\mathbb{C}^2)\right)-\cosh\left(\frac
e2\right)\mathrm{ch}\left(B_j(T_\mathbb{C}B+N_\mathbb{C},
N_\mathbb{C})\right)
}{2\sinh\left(\frac{e}{2}\right)}\right\}^{(8k+6)}, $$ $0\leq
j\leq r.$ Actually, $h_r$ is just
$$\left\{\widehat{A}(TB,\nabla^{TB}){{\rm ch}\left(
b_r(T_{\mathbb C}B+N_{\mathbb C},{\mathbb C }^2)\right)
-\cosh\left({e\over 2}\right){\rm ch}\left( b_r(T_{\mathbb
C}B+N_{\mathbb C},N_{\mathbb C })\right)\over 2\sinh\left({e\over
2}\right)}\right\}^{(8k+6)}.$$ \end{theorem}

\subsection {A general type of cancellation formulas for even
dimensional Riemannian manifolds} In this subsection, let's
continue to discuss a general type of cancellation formulas. Let
$B$ be an $2d$-dimensional smooth oriented Riemannian manifold,
$m$ be a non-negative integer and $N$ be a complex line bundle on
$B$. Define \be \Theta_1'(T_\mathbb{C}B, m,
N_\mathbb{C})=\bigotimes_{n=1}^\infty
S_{q^n}(\widetilde{T_\mathbb{C}B}-m\widetilde{N_\mathbb{C}})
\otimes\bigotimes_{m=1}^\infty\Lambda_{q^m}
(\widetilde{T_\mathbb{C}B}-m\widetilde{N_\mathbb{C}}), \ee \be
\Theta_2'(T_\mathbb{C}B, m, N_\mathbb{C})=\bigotimes_{n=1}^\infty
S_{q^n}(\widetilde{T_\mathbb{C}B}-m\widetilde{N_\mathbb{C}})
\otimes\bigotimes_{s=1}^\infty\Lambda_ {-q^{s-\frac 12}}
(\widetilde{T_\mathbb{C}B}-m\widetilde{N_\mathbb{C}}), \ee and
assume they admit Fourier expansion in the following:

\be \Theta_1'(T_{\mathbb C}B, m ,N_{\mathbb C})=A_0'(T_{\mathbb
C}B, m, N_{\mathbb C})+ A_1'(T_{\mathbb C}B, m, N_{\mathbb
C})q^{1/2}+\cdots, \ee \be\Theta_2'(T_{\mathbb C}B, m, N_{\mathbb
C})=B_0'(T_{\mathbb C}B, m, N_{\mathbb C}) +B_1'(T_{\mathbb C}B,
m, N_{\mathbb C})q^{1/2}+\cdots,
 \ee

Set \be Q_1'(\tau)=\widehat{L}(TB,\nabla^{TB})
\frac{\sinh^{2n+\frac{1-(-1)^d}{2}}\left(\frac
e2\right)}{\cosh^{2n+\frac{1-(-1)^d}{2}}\left(\frac
e2\right)}\mathrm{ch}(\Theta_1'(T_\mathbb{C}B,
2n+\frac{1-(-1)^d}{2}, N_\mathbb{C})),\ee

\be
Q_2'(\tau)=\widehat{A}(TB,\nabla^{TB})\sinh^{2n+\frac{1-(-1)^d}{2}}\left(\frac
e2\right)\mathrm{ch}(\Theta_2'(T_\mathbb{C}B,
2n+\frac{1-(-1)^d}{2}, N_\mathbb{C})),\ee where $n$ is a
nonnegative integer and satisfy
$d-\left(2n+\frac{1-(-1)^d}{2}\right)>0$.

By similar computations as those in the proof of Theorem 3.1, we
have \be
\begin{split} Q_1'(\tau)=&2^{d}\left(\prod_{j=1}^{d}x_j\frac{\theta'(0,\tau)}{\theta(x_j,\tau)}
          \frac{\theta_1(x_j,\tau)}{\theta_1(0,\tau)}\right)\\
     & \cdot{\left(\frac{\theta(u,\tau)}{\theta'(0,\tau)}\right)}^{2n+\frac{1-(-1)^d}{2}}
{\left(\frac{\theta_1(0,\tau)}{\theta_1(u,\tau)}\right)}^{2n+\frac{1-(-1)^d}{2}},
\end{split}\ee

\be \begin{split}
Q_2'(\tau)=&\left(\prod_{j=1}^{d}x_j\frac{\theta'(0,\tau)}{\theta(x_j,\tau)}
                                             \frac{\theta_2(x_j,\tau)}{\theta_2(0,\tau)}\right)\\
                  & \cdot{\left(\frac{\theta(u,\tau)}{\theta'(0,\tau)}\right)}^{2n+\frac{1-(-1)^d}{2}}
                  {\left(\frac{\theta_2(0,\tau)}{\theta_2(u,\tau)}\right)}^{2n+\frac{1-(-1)^d}{2}}.
\end{split}\ee
Let $P_1'(\tau)=\{Q_1'(\tau)\}^{(2d)},
P_2'(\tau)=\{Q'_2(\tau)\}^{(2d)}$ be the $(2d)$-components of
$Q_1'(\tau), Q'_2(\tau)$ respectively. $P_1'(\tau)$ is a modular
form of weight $d-\left(2n+\frac{1-(-1)^d}{2}\right)$ over
$\Gamma_0(2)$ and $P_2'(\tau)$ is a modular form of weight
$d-\left(2n+\frac{1-(-1)^d}{2}\right)$ over $\Gamma^0(2)$. Playing
the same game as in the proof of Theorem 3.1, we obtain
\begin{theorem} The following identity holds,

\be \begin{split}
\frac{1}{2^{\frac{3d}{2}-\frac{1-(-1)^d}{4}-n}}\left\{\widehat{L}(TB,\nabla^{TB})\frac{\sinh^{2n+\frac{1-(-1)^d}{2}}\left(\frac
e2\right)}{\cosh^{2n+\frac{1-(-1)^d}{2}}(\frac
e2)}\right\}^{(2d)}\\
=\sum_{r=0}^m2^{-6r}\left\{d_r(B,2n+\frac{1-(-1)^d}{2},N)\sinh^{2n+\frac{1-(-1)^d}{2}}\left(\frac{e}{2}\right)\right\}^{(2d)}\,
,
\end{split}\ee
where each $d_r(B,2n+\frac{1-(-1)^d}{2},N), 0\leq r\leq k$, is a
finite and canonical linear combination of characteristic forms
$\widehat{A}(TB,\nabla^{TB}) \mathrm{ch}\left(B_i'(T_\mathbb{C}B,
2n+\frac{1-(-1)^d}{2}, N_\mathbb{C})\right),\ 0\leq i\leq r $ and
$m={\left[\frac{d-2n-\frac{1-(-1)^d}{2}}{4}\right]}$.\end{theorem}

\begin{remark} The condition
$d-\left(2n+\frac{1-(-1)^d}{2}\right)>0$ is put to make Theorem
3.3 nontrivial. If $d-\left(2n+\frac{1-(-1)^d}{2}\right)=0$, then
both sides are
$\frac{e^{2d}}{2^{\frac{3d}{2}-\frac{1-(-1)^d}{4}-n}}$. If
$d-\left(2n+\frac{1-(-1)^d}{2}\right)<0$, then both sides of
Theorem 3.3 are zeros since $\sinh$ is an odd function and the
degrees of the top components of both sides are greater than $2d$.
\end{remark}

\begin{remark} When $n=0$ and $\frac{1-(-1)^d}{2}=1$, i.e. $d=4a+1$
or $4a+3$, the integral of the left hand side of Theorem 3.3
against the fundamental class of $B$ is, up to a constant, the
signature of a submanifold of $B$ which is the smooth zero locus
of a generic section of the bundle $N$. Thus when $B$ is a
$spin^c$ manifold and $N$ is dual to $w_2(B)$, Theorem 3.3 shows
that the signature of the smooth submanifolds of B dual to $e(N)$
can be given by indexes of Dirac operators on the $spin^c$
manifold B.

\end{remark}

Putting $d=4k+2$ and $n=0$ in Theorem 3.3, we get
\begin{corollary} The following cancellation formula holds,
$$\Big\{\widehat{L}(TM,\nabla^{TM})\Big\}^{(8k+4)}=8\sum_{j=0}^k
2^{6k-6j}\Big\{\widehat{A}(TM,\nabla^{TM})
\mathrm{ch}b_j\Big\}^{(8k+4)},$$ where the $b_j$'s are elements in
$KO(M)\otimes {\mathbb C}$.
\end{corollary}
This is just Liu's original cancellation formula [10]. So Theorem
3.3 is a generalization of Liu's cancellation formula to all even
dimensional oriented Riemannian manifolds with a complex line
bundle involved. In particular, when $d=6, n=0$, we get the
Alvarez-Gaum\'e-Witten miraculous cancellation formula (1.1).

Looking at Theorem 3.3, let's get some interesting cancellation
formulas for special $d$ and $n$. Putting $d=6$ and $n=1$, i.e.
for 12-dimensional manifold $M$, we have
\begin{corollary} The following formula holds, \be
\left\{\widehat{L}(TM,
\nabla^{TM})\frac{\sinh^2{\left(\frac{c}{2}\right)}}{\cosh^2{\left(\frac{c}{2}\right)}}\right\}^{(12)}=\left\{
\left[ -4\widehat{A}(TM, \nabla^{TM}){\rm ch} (
T_{\mathbbm{C}}M,\nabla^{T_{\mathbbm{C} }M})\right.\right.\ee
$$+112\widehat{A}(TM,
\nabla^{TM})+8\left.\left.\widehat{A}(TM,
\nabla^{TM})\left(e^c+e^{-c}-2\right)\right] \sinh^2\left({c\over
2}\right)\right\}^{(12)}.$$
\end{corollary}

Putting $d=6$ and $n=2$, i.e. for 12-dimensional manifold $M$, we
have
\begin{corollary} The following formula holds, \be
\left\{\widehat{L}(TM,
\nabla^{TM})\frac{\sinh^4{\left(\frac{c}{2}\right)}}{\cosh^4{\left(\frac{c}{2}\right)}}\right\}^{(12)}=-128\left\{\widehat{A}(TM,
\nabla^{TM})\sinh^4{\left(\frac{c}{2}\right)}\right\}^{(12)}.\ee
\end{corollary}

Corollary 3.2 and 3.3 are both analogous to the
Alvarez-Gaum\'e-Witten original miraculous cancellation formula
(1.1) with a complex line bundle involved.

Putting $d=5$ and $n=0$, i.e. for 10-dimensional manifold $M$, we
have
\begin{corollary} The following formula holds, \be
\left\{\widehat{L}(TM,
\nabla^{TM})\frac{\sinh{\left(\frac{c}{2}\right)}}{\cosh{\left(\frac{c}{2}\right)}}\right\}^{(10)}=\left\{
\left[ -2\widehat{A}(TM, \nabla^{TM}){\rm ch} (
T_{\mathbbm{C}}M,\nabla^{T_{\mathbbm{C} }M})\right.\right.\ee
$$+52\widehat{A}(TM,
\nabla^{TM})+2\left.\left.\widehat{A}(TM,
\nabla^{TM})\left(e^c+e^{-c}-2\right)\right] \sinh\left({c\over
2}\right)\right\}^{(10)}.$$
\end{corollary}

Putting $d=5$ and $n=1$, i.e. for 10-dimensional manifold $M$, we
have
\begin{corollary} The following formula holds, \be
\left\{\widehat{L}(TM,
\nabla^{TM})\frac{\sinh^3{\left(\frac{c}{2}\right)}}{\cosh^3{\left(\frac{c}{2}\right)}}\right\}^{(10)}=-64\left\{\widehat{A}(TM,
\nabla^{TM})\sinh^3{\left(\frac{c}{2}\right)}\right\}^{(10)}.\ee
\end{corollary}

$\ $

\section {Relations among the cancellation formulas}
In this section, by applying the method of integration along the
fibre, we show some relations on the level of characteristic
numbers among cancellation formulas obtained in [10], [7] and \S 3
of this article.

First let's get (3.1) from Theorem 2.1 by integration along the
fibre. Let $u \in H^2_{cv}(N)$, the second compact vertical
supports cohmology of $N$, be the Thom class of the bundle $(N,
\pi, B)$ with fibre $L$. By the Thom isomorphism theorem, we have
the following identity of cohomology classes in $H^*(B)$,
 $$\left[\int_Lu^{2i}\right]=\left[e^{2i-1}\right],\, i=1, 2,\cdots.$$

By integration along the fibre, on the one hand, we have \be
\begin{split}&\int_{N}\widehat{L}(TN)\left(1-\frac{1}{\cosh^{2}\left(\frac{u}{2}\right)}\right)\\
=&\int_B\widehat{L}(TB)\int_{L}\frac{u}{\sinh\left(\frac{u}{2}\right)}\cosh\left(\frac{u}{2}\right)
\frac{\sinh^2\left(\frac{u}{2}\right)}{\cosh^2\left(\frac{u}{2}\right)}\\
=&\int_B\widehat{L}(TB)\frac{\sinh\left(\frac{e}{2}\right)}
{\cosh\left(\frac{e}{2}\right)}\end{split}\ee

\noindent and on the other hand,  $\forall 0\leq r \leq k,$ we
have \be \begin{split}& \int_{N}\widehat{A}(TN)
\left(\mathrm{ch}\left(b_r(T_\mathbb{C}N,\mathbb{C}^2)\right)
-\cosh\left(\frac{u}{2}\right)\mathrm{ch}\left(b_r(T_\mathbb{C}N,\xi_{\mathbb{C}})\right)\right)\\
=&\int_B \widehat{A}(TB)\int_{L}{{\sinh\left({u\over
2}\right)}\over {\cosh\left({u\over 2}\right)} }
\left\{\mathrm{ch}\left(b_r\left(T_\mathbb{C}N,\mathbb{C}^2\right)\right)
-\cosh\left(\frac{u}{2}\right)\mathrm{ch}\left(b_r\left(T_\mathbb{C}N,\xi_{\mathbb{C}}
\right)\right)\right\}\\
=&\int_B
\widehat{A}(TB)\frac{\mathrm{ch}\left(b_r(T_\mathbb{C}B+N_\mathbb{C},
\mathbb{C}^2)\right)-\cosh\left(\frac{e}{2}\right)\mathrm{ch}\left(b_r(T_\mathbb{C}B+N_\mathbb{C},
N_\mathbb{C})\right)}{2\sinh\left(\frac{e}{2}\right)}.
\end{split}\ee
Then by Theorem 2.1, Corollaries 2.1, Corollary 2.2 and (4.1),
(4.2), we obtain
$$
 { 1 \over 8}\int_B \widehat{L}(TB){{\sinh\left({e\over
2}\right)}\over {\cosh\left({e\over 2}\right)} }$$
$$=\sum_{r=0}^k 2^{6k-6r}\int_B\widehat{A}(TB){{\rm
ch}\left( b_r(T_{\mathbb C}B+N_{\mathbb C},{\mathbb C }^2)\right)
-\cosh\left({e\over 2}\right){\rm ch}\left( b_r(T_{\mathbb
C}B+N_{\mathbb C},N_{\mathbb C })\right)\over 2\sinh\left({e\over
2}\right) }.$$

\begin{remark} one can also obtain the characteristic number
version of Theorem 3.2 from an $8k$-analogue of Theorem 2.1 stated
in [7, Theorem A.1] by integration along the fibre as
above.\end{remark} Next let's say something about the relations
among formulas in Theorem 3.3 for different $d$ and $n$.

Let $u \in H^2_{cv}(N)$, the second compact vertical supports
cohmology of $N$, be the Thom class of the bundle $(N, \pi, B)$
with fibre $L$. By the Thom isomorphism theorem, we have the
following identity of cohomology classes in $H^*(B)$,
 $$\left[\int_Lu^{2i}\right]=\left[e^{2i-1}\right],\, i=1, 2,\cdots.$$

Let $M$ be an $(8k+4)$-dimensional closed oriented Riemannian
manifold and $(E, \pi, M)$ be a complex line bundle on $M$ with
fibre $L$. Let $u \in H^2_{cv}(E)$ be the Thom class of this
bundle. Then from Theorem 3.3 in the case $d=4k+3$ and $n=0$, one
has

$$
\left\{\widehat{L}(TE,\nabla^{TE})\frac{\sinh\left(\frac
u2\right)}{\cosh(\frac u2)}\right\}^{(8k+6)}
=16\sum_{r=0}^{k}2^{6k-6r}\left\{d_r(E, 1,
\pi^*E)\sinh\left(\frac{u}{2}\right)\right\}^{(8k+6)}\, ,
$$
where each $d_r(E, 1, \pi^*E),0\leq r\leq k$, is a finite and
canonical linear combination of characteristic forms
$\widehat{A}(TE,\nabla^{TE}) \mathrm{ch}\left(B_i'(T_\mathbb{C}E,
1, {(\pi^*E)}_\mathbb{C})\right),\ 0\leq i\leq r$. Performing
integration along the fibre, we have
$$\int_E \widehat{L}(TE)\frac{\sinh\left(\frac
u2\right)}{\cosh(\frac u2)}=\int_M \widehat{L}(TM)\int_L
\frac{u}{\tanh \left({u \over 2}\right)}\frac{\sinh\left(\frac
u2\right)}{\cosh(\frac u2)}=\int_M \widehat{L}(TM)$$ and
$$\int_E \widehat{A}(TE)
\mathrm{ch}\left((B_i'(T_\mathbb{C}E, 1,
{(\pi^*E)}_\mathbb{C}))\right)\sinh\left(\frac u2\right)$$
$$=\int_M \widehat{A}(TM)\int_L \frac{{u\over 2}}{\sinh\left(\frac
u2\right)}\mathrm{ch}\left((B_i'(T_\mathbb{C}E, 1,
{(\pi^*E)}_\mathbb{C})\right)\sinh\left(\frac u2\right)$$
$$={1\over 2}\int_M
\widehat{A}(TM)\mathrm{ch}\left((B_i'(T_\mathbb{C}M, 0,
E_\mathbb{C})\right).
$$
Hence we obtain
$$ \int_M\widehat{L}(TM)=8\sum_{j=0}^k
2^{6k-6j}\int_M \widehat{A}(TM) \mathrm{ch}d_j(M, 0, E),
$$
which is just the characteristic number version of the case of
$d=4k+2$ and $n=0$ in Theorem 3.3.

More generally, with the same pattern, we can apply integration
along the fibre to get the formula in the case of $(d,n)$ from the
formula in the case of $(d+1, n+\frac{1-(-1)^d}{2})$ on the level
of characteristic numbers. This phenomena looks very interesting
since it beautifully relates different cancellation formulas in
Theorem 3.3 which are all products of modular invariance.

\section {Acknowledgments} The authors would like to thank
Professor Weiping Zhang for his very helpful suggestions and
encouragement. The results of this paper were once reported in a
seminar held in the Department of Mathematics of UCSB. The authors
thank Professor Xianzhe Dai for his interest and encouragement. We
are also grateful to Professor Huitao Feng for his caring in the
writing of this paper.

\bibliographystyle{amsplain}

\begin{thebibliography}{10}


\bibitem {A} L. Alvarez-Gaum\'e and E. Witten, Gravitational
anomalies. {\it Nucl. Phys.} B234 (1983), 269-330.

\bibitem {A} M. F. Atiyah, $K-theory$. Benjamin, New York, 1967.

\bibitem {A} M. F. Atiyah and F. Hirzebruch, Riemann-Roch
theorems for differentiable manifolds. {\it Bull. Amer. Math.
Soc.} 65 (1959), 276-281.

\bibitem {C} K. Chandrasekharan, {\it Elliptic Functions}. Springer-Verlag,
1985.

\bibitem {F} S. M. Finashin, A Pin$^-$-cobordism invariant and a
generalization of Rokhlin signature congruence. {\it Leningrad
Math. J.} 2 (1991), 917-924.

\bibitem {H} F. Han and W. Zhang, Spin$^{c}$-manifold and elliptic genera. {\it
C. R. Acad. Sci. Paris, S$\acute{e}$rie I.} 336 (2003), 1011-1014.
\bibitem {H} F. Han and W. Zhang, Modular invariance, characteristic numbers
and $\eta$ invariants. {\it Journal of Differential Geometry.} 67
(2004), 257-288.

\bibitem {H} F. Hirzebruch, {\it Topological Methods in
Algebraic Geometry.} Springer-Verlag, 1966.

\bibitem {L} P. S. Landweber, Elliptic cohomology and modular forms. in {\it
Elliptic Curves and Modular Forms in Algebraic Topology, } p.
55-68. Ed. P. S. Landweber. Lecture Notes in Mathematics Vol.
1326, Springer-Verlag (1988).

\bibitem {L} K. Liu, Modular invariance and characteristic
numbers. {\it Commun. Math. Phys}. 174 (1995), 29-42.

\bibitem {L} K. Liu and W. Zhang, Elliptic genus and
$\eta$-invariants. {\it Inter. Math. Res. Notices} No. 8 (1994),
319-328.

\bibitem {O} S. Ochanine, Signature modulo 16, invariants de
Kervaire g\'eneralis\'e et nombre caract\'eristiques dans la
$K$-th\'eorie reelle. {\it M\'emoire Soc. Math. France}, Tom. 109
(1987), 1-141.

\bibitem {Z} W. Zhang, Spin$^c$-manifolds and Rokhlin
congruences.
 {\it   C. R. Acad. Sci. Paris,
S\'erie I}, 317 (1993), 689-692.

\bibitem {Z} W. Zhang, Circle bundles, adiabatic limits of
$\eta$ invariants and Rokhlin cogruences. {\it Ann. Inst. Fourier
} 44 (1994), 249-270.

\bibitem {Z} W. Zhang, {\it Lectures on Chern-Weil Theory and
Witten Deformations.} Nankai Tracts in Mathematics Vol. 4, World
Scientific, Singapore, 2001.
$$\  $$

\end{thebibliography}

\end{document}